\documentclass[12pt]{article}
\usepackage{amssymb}
\usepackage[cp1250]{inputenc}
\usepackage[T1]{fontenc}
\usepackage{hyperref}
\def\N{{\mathbb N}}
\def\Q{{\mathbb Q}}
\def\b{{\rm {\bf b}}}
\def\d{{\rm {\bf d}}}
\def\P{{\N^\N}}

\title{\large \bf
A remark on  Fremlin--Miller theorem concerning the \nolinebreak Menger property
and Michael concentrated sets}

\author
{\large \bf J.~Chaber, R.~Pol }

\date{10.10.2002}

\begin{document}
\maketitle
\begin{abstract}
The theorem we prove is a slight strengthening of some results by
Just, Miller, Scheepers and Szeptycki \cite{JMSS}. We use the
Michael technique instead of the combinatorial approach in the
literature.
\end{abstract}

We consider only separable metrizable spaces. A space $M$ has the
Menger property if for any sequence ${\mathcal U}_1,{\mathcal U}_2,\ldots$
of open covers of $M$ there are finite collections ${\mathcal
V}_i\subset {\mathcal U}_i$ with $\bigcup_i{\mathcal V}_i$ covering $M$.
Let $\d$ be the minimal cardinality of a dominating set in $\P$
(we recall the definition below).

Fremlin and Miller \cite{FM} proved that any space of cardinality
less than $\d$ has the Menger property. In particular, upon an
identification of $\P$ with the irrationals, if $S$ is a scale in
$\P$ of cardinality $\d$, and $\Q$ are the rationals, the union
$S\cup\Q$ has the Menger property (notice also that $S$ is
concentrated about $\Q$). In effect, no matter where $\d$ is
situated on the alephs scale, the Fremlin--Miller theorem yields a
non-$\sigma$-compact space with the Menger property.

We shall add to this important result an observation that the
Michael technique of constructing concentrated sets in product
spaces, introduced in \cite{Mi}, allows one to use the
Fremlin--Miller theorem to the following effect.

\medskip
{\bf Theorem.} {\it In every completely metrizable, separable
space $X$ without any compact set with nonempty interior there is
an $M\subset X$ whose all finite powers have the Menger property,
but $M$ is not contained in any $\sigma$-compact subset of $X$.}
\medskip

Some closely related results were proved by Just, Miller,
Scheepers and Szeptycki with the use of a different, combinatorial
approach, cf.\ Theorems 3.9 and 5.1 in \cite{JMSS}.

Before passing to a proof of the theorem, let us recall some facts
about the order $\leq^*$ in $\P$, where $s\leq^* t$ means that
$s(n)\leq t(n)$ for all but finitely many $n$. The minimal
cardinality of a set in $\P$ unbounded with respect to $\leq^*$ is
denoted by $\b$, and $\d$ is the minimal cardinality of a set
dominating in $\P$, i.e., a set $D\subset\P$ such that for any
$s\in \P$ there is a $t\in D$ with $s\leq^* t$, cf.\ \cite{vD}.
Always $\b\leq\d$, and the equality $\b=\d$ yields a dominating
set of cardinality $\d$ well-ordered by $\leq^*$, i.e., a scale in
$\P$.

\medskip
{\bf Proof of Theorem.} To begin with, let us recall that each
$t\in \P$ determines a $\sigma$-compact set
$F(t)=\{s\in\P:s\leq^*t\}$ and each $\sigma$-compact set in $\P$ is
contained in some $F(t)$.

Let $X$ satisfy the assumptions and let us fix a countable set $Q$
dense in $X$.

Let \begin{equation} f:\P\rightarrow X \end{equation} be a perfect
surjection.

Assume $\b<\d$, let $B\subset\P$ be an unbounded set of
cardinality $\b$ and let $M=f(B)$. Since the cardinality of $M$ is
less than $\d$, any finite power of $M$ has the Menger property by
Fremlin--Miller theorem and $M$ is not contained in any
$\sigma$-compact set in $X$.

From now on we assume $\b=\d$. Then we have a scale
$\{t_\alpha:\alpha<\d\}$ in $\P$, $t_\alpha \leq^* t_\beta$ for
$\alpha<\beta$. Let, cf.\ (1),
 \begin{equation}
F_\alpha=f(F(t_\alpha)).
 \end{equation}
Each $\sigma$-compact set in $X$ is contained in some $F_\alpha$.
Let $Z$ be a compactification of $X$ and let $$g_m:\P\rightarrow
Z^m\setminus Q^m$$ be perfect surjections, $H_\alpha=X^m\setminus
g_m(F(t_\alpha))$. Then $H_\alpha\supset Q^m$ and, for any
collection ${\mathcal H}$ of $G_\delta$ sets in $X^m$ containing $Q^m$,
if $|{\mathcal H}|<\d$ then there is $\alpha$ with
$H_\alpha\subset\bigcap{\mathcal H}$.

In effect, changing the enumeration of the sets $H_\alpha$, we get
$G_\delta$-sets $G_\xi$ with
\begin{equation}
Q^{m(\xi)}\subset G_\xi\subset X^{m(\xi)},\quad \xi<\d,
\end{equation}
such that for any collection ${\mathcal H}$ of $G_\delta$-sets in $X^m$
containing $Q^m$,
\begin{equation}
 \mbox{ if
}|{\mathcal H}|<\d \mbox{ then }G_\xi\subset\bigcap{\mathcal H}
 \mbox{ for some }\xi \mbox{ with } m(\xi)=m.
\end{equation}

We shall define
\begin{equation}
B=\{x_\alpha:\alpha<\d\}\subset X\setminus Q,\quad \mbox{and}\quad
M=Q\cup B.
\end{equation}
such that the following concentration property holds true (cf.\
\cite[Lemma 2.2.2]{PPR}): for any $m$ and each collection ${\mathcal
H}$ of $G_\delta$ sets in $X^m$ containing $Q^m$, if $|{\mathcal
H}|<\d$ then $M^m\setminus \bigcap{\mathcal H}$ is contained in the
union of less than $\d$ hyperplanes obtained by fixing one of the
coordinates in $M^m$.

A straightforward induction on $m$ shows that for such $M$ every
product $L\times M^m$, where $|L|<\d$, has the Menger property
(for $m=0$ use the Fremlin--Miller theorem). Therefore it is
enough to define the set $B$ in (5), ensuring also that $B$ is not
contained in any $\sigma$-compact set in $X$.

To this end, we shall modify slightly a reasoning in \cite[Comment
9.5]{PPR}, which is a variation of the celebrated Michael technique
\cite{Mi}.

We shall start with some notation. Let $S\subset\{1,2,...,m\}$,
$u\in X^S$, $v\in X^{\{1,2,...,m\}\setminus S}$. Then $(u,v)$ is
the point $w \in X^m$ with $w(i)=u(i)$ for $i\in S$ and $w(i)=v(i)$
for $i\notin S$. Further, given a $G_\delta$-set $G$ in $X^m$
containing $Q^m$, we let $$G(S)=\{u\in X^S: (u,v)\in G \mbox{ for
all } v\in Q^{\{1,2,...,m\}\setminus S}\},$$ and, for $u\in G(S)$,
$j\in{\{1,2,...,m\}\setminus S}$, we let $$G(u,j)=\{x\in X:
(u,x)\in G(S\cup\{j\})\}.$$ Then $G(S)$, $G(u,j)$ are
$G_\delta$-sets with
\begin{equation}
Q^S\subset G(S)\subset X^S,\quad  Q\subset G(u,j)\subset X
\end{equation}

We shall pick distinct points $x_\alpha$ of the set $B$ in (5) by
transfinite induction, such that, with $\Delta(S)=\{u\in X^S:
\mbox{ at least two coordinates of } u \mbox{ are equal}\},$ we
shall have for any $\alpha<\d$, $\xi<\alpha$, $S\subset
\{1,2,...,m(\xi)\}$,
$$
\{x_\eta:\xi<\eta\leq\alpha\}^S\setminus\Delta(S)\subset G_\xi(S)
$$

We begin with an arbitrary $x_0\notin Q$. At each stage
$\alpha>0$, $x_\alpha$ must be chosen from the intersection of a
collection ${\mathcal H}$ consisting of sets of the form $G_\xi(u,j)$,
$\xi<\alpha$, $u\in G_\xi(S)\cap \{x_\eta:\xi<\eta<\alpha\}^S
\setminus\Delta(S)$, $S\cup\{j\}\subset\{1,2,...,m(\xi)\}$, and
$j\notin S$. But these are $G_\delta$-sets in $X$ containing $Q$,
cf.\ (6), (3), and $|{\mathcal H}|<\d$, hence by (4), the intersection
$\bigcap{\mathcal H}$ contains a $G_\delta$-set in $X$ containing $Q$.
This shows that a choice of $x_\alpha$ is always possible and,
moreover, one can pick $x_\alpha$ outside of the set $F_\alpha$,
meager in $X$, cf.\ (2), which will also guarantee that the
resulting set $B$ is not contained in any $\sigma$-compact set in
$X$. A verification that the space $M=Q\cup B$ has the
concentration property described below (5) is similar to a simple
reasoning in \cite[Comment 9.5]{PPR}.


\begin{thebibliography}{ChGP??}

\bibitem[vD]{vD} E.K. van Douwen, {\it The Integers and Topology\/},
         Handbook of Set-Theoretic Topology, K. Kunen and J.E.Vaughan,
         editors, Elsevier Science Publishers, 1984, 111--167.

\bibitem[JMSS]{JMSS} W. Just, A.W. Miller, M. Scheepers, P.J. Szeptycki,
        {\it The combinatorics of open covers II\/}, Topology and Appl.
        73(1996), 241--266.

\bibitem[Mi]{Mi} E.A. Michael, {\it Paracomactness and the Lindel\"{o}f
        property in finite and countable cartesian products\/},
        Compositio Math. 23(1971), 199--244.

\bibitem[MF]{FM} A.W. Miller, D.H. Fremlin {\it On some properties of Hurewicz,
    Menger, and Rothberger\/}, Fund. Math. 129(1988), 17--33.

\bibitem[PPR]{PPR} E. Pol, R. Pol, M. Re\'{n}ska, {\it On countable-dimensional
        spaces with the Menger property, rational dimension and a question
        of S.D. Iliadis\/}, Mh. Math. 128(1999), 331--348.

\end{thebibliography}
\end{document}